\documentclass{amsart} 
\usepackage{bbm,longtable,graphicx}
\title{Affine extension of noncrystallographic Coxeter groups 
and quasicrystals}

\author{Ji\v r\'\i\  Patera}
\thanks{Work supported in part by the Natural Sciences and
Engineering Research Council of Canada}
\address{Centre de Recherches Math\'ematiques\\ 
  Universit\'e de Montr\'eal, Montr\'eal, Qu\'ebec, Canada}
\email{patera@crm.umontreal.ca}

\author{Reidun Twarock}
\thanks{Marie Curie Fellow} 
\address{Department of Mathematics\\ 
University of York, YO10 5DD, England}
\email{rt11@york.ac.uk}

\subjclass{02.20.s, 61.44}

\keywords{noncrystallographic Coxeter groups, affine
extension, quasicrystals}

\newtheorem{lem}{\bf Lemma}[section] 
\newtheorem{de}[lem]{\bf Definition} 
\newtheorem{thm}[lem]{Theorem} 
 
\newtheorem{prop}[lem]{Proposition} 
\newtheorem{coro}[lem]{Corollary} 

\def\t{\tau}

\def\R{\mathbb{R}}
\def\C{\mathbb{C}} 
\def\Z{\mathbb{Z}} 
\def\zt{\mathbb{Z}[\tau]} 
\def\N{\mathbb{N}} 
\def\Q{\mathbb{Q}} 
\def\integer{\mathbb{Z}} 
\def\half{\frac{1}{2}} 
\def\aff{\mathop{\rm aff}}
 
\begin{document} 
\begin{abstract} Unique affine extensions $H^{\aff}_2$,
$H^{\aff}_3$ and $H^{\aff}_4$ are determined for the
noncrystallographic Coxeter groups $H_2$, $H_3$  and
$H_4$. They are used for the construction of new mathematical
models for quasicrystal fragments with 10-fold symmetry. The case
of $H^{\aff}_2$ corresponding to planar point sets is
discussed in detail. In contrast to the cut-and-project
scheme we obtain by construction finite point sets,
which grow with a model specific growth parameter. 
\end{abstract}

\maketitle


\section{Introduction} 

In contrast to the well known Weyl groups of affine
Kac-Moody algebras, affine extensions of finite noncrystallographic
Coxeter groups have apparently not been studied, although they too are a
natural part of the general theory of Coxeter groups of infinite order. 
The goal of this article is to describe such extensions for the Coxeter
groups $H_2$, $H_3$, and $H_4$ of order $10$, $120$, and $14\,440$, 
respectively. Since our motivation for this study comes from the
theory of quasicrystals, we illustrate the exploitation of
such groups on problems related to quasicrystal generation/growth, but 
we expect applications also to other areas such as e.g. fullerenes 
\cite{onion} as we comment below. 

There is a far going parallel between the finite noncrystallogaphic and the
crystallographic Coxeter groups, the latter being the Weyl groups
associated with simple Lie algebras. In both cases, the groups are
 generated by reflections and have a unique affine extension. The main
difference between the two consists in the fact that the crystallographic
affine groups generate an entire root lattice starting from any root or
from the origin, whereas similar applications of a noncrystallographic group
to the origin or to any root of  $H_2$, $H_3$, and $H_4$ would generate
a point set which densely covers the whole space. 

There is an important property which makes the $H_k$-cases richer
than the crystallographic ones: There exists a root map, that is a 
mapping transforming root systems into root systems, which is not from 
the Coxeter group, and which acts as a nontrivial
transformation of $H_k$-roots. It is 
the mapping  called the star map in \cite{CMP} which, for example, provides the one-to-one correspondence 
between quasicrystal points and the points in the corresponding acceptance 
window.

Since the discovery of quasicrystals in physics \cite{Shechtman:1984},
mathematical models describing these aperiodic structures have been
proposed. Perhaps the best established is the cut-and-project approach
for the construction of point sets  modelling quasicrystals \cite{Duneau}. Through the
years a number of variants of the method have been developed (see for
example \cite{Janot:1994} and references therein). Our considerations are
based on an algebraic way of construction 
\cite{CMP,P,Moody:1997b,Moody:2000},  in which the uniformity of the
procedure for different dimensions allows to consider these models 
simultaneously subject only to a variation of the starting data.
Properties of the cut-and-project point sets are now understood in
great details particularly in one dimension
\cite{Patera:1998a,Patera:1998b,Patera:1999,Patera:1998c}. A key
constituent in models related to point sets with 10-fold rotational
symmetry are the noncrystallographic Coxeter groups  $H_2$, $H_3$ and
$H_4$, leading to models in 2, 3 and 4 dimensions, respectively.  They
exploit the fact that
$H_2$, $H_3$ and $H_4$ point sets are projections from crystallographic lattices
of types
$A_4$,
$D_6$ or
$E_8$,  respectively \cite{P,Moody:1993}. 

Coxeter groups \cite{Cox,Humphreys:1992} are discrete groups generated
by reflections. A special class among them are the Weyl groups (or
crystallographic Coxeter groups) which are the finite symmetry groups of
root and weight lattices in the theory of semisimple Lie algebras/groups
and their representations. Affine extensions of the Weyl groups are also
generated by reflections. They are of infinite order and are known to
underlie similar symmetries of the affine Kac-Moody algebras
\cite{Kac:1985,Kass:1990}. Finite noncrystallographic Coxeter groups
(which are not products of several smaller ones) can be easily enumerated.
Those generated by more than two reflections are two:
$H_3$ generated by three reflections, and $H_4$ generated by four reflections.
Coxeter groups generated by  two reflections are infinitely many: they
are the symmetry groups of regular polygons with any number of vertices
but
2, 3, 4, and 6. (The latter ones are of crystallographic type). In
addition to $H_3$ and $H_4$, it is natural to
consider in this article also the lowest of the 2-reflection groups,
called
here $H_2$, the symmetry group of regular pentagons and decagons. In more
familiar physics terminology $H_2$ is the dihedral group of order 10,
while $H_3$ is the icosahedral group of order 120. A  description of the
three groups
suitable for our problem can be found in Section~2.

The group $H_4$, which is of order $14\,400$, does not have a standard
name
in physics, nevertheless on a few occasions it appeared in the physics
literature either in the context of the physics of amorphous solids
\cite{Di-1,Di-2,Di-3,Di-4,Di-5,Frad}, biophysics \cite{Bul'en}, 
quasicrystals \cite{Sloane,Moody:1993}, or general mathematical
physics \cite{Lam,CKPS}. 
The group $H_4$ contains all point groups familiar in 3-dimensional
crystallography, besides the inclusions $H_4\supset H_3\supset H_2$.
Moreover there is the remarkable relation (see \cite{Sch,Moody:1993,
CMP,P}
and references therein) of $H_4$ to the largest exceptional simple Lie
group $E_8$ encountered in particle physics \cite{Gross}. Therefore it is
possible that $H_4$ and/or its affinisation $H^{\aff}_4$ will  play a
basic role in physics in not too distant future.

For further information about the non-extended groups, see e.g.
\cite{CMP,P,Humphreys:1992}.  The diagrams representing our affine
extended groups correspond to graphs related to regular polytopes
\cite{Cox}. We remark that different generalizations of
finite Coxeter groups and related diagrams  appear also in
\cite{Zuber:1995,Zuber:1997}.

 Affine extensions $H^{\aff}_2$, $H^{\aff}_3$, and $H^{\aff}_4$ of the
groups $H_2$, $H_3$, and $H_4$, unlike the affine extensions of the Weyl
groups, have apparently been considered neither in the physics nor in
the mathematics literature before. In this article we first describe the
affine extensions of the three groups, pointing out particularly their
uniqueness and the close analogy to the affinisation of the Weyl
groups. 

In Section~2 we  recall  the way in which root lattices are
constructed in Lie theory based on affine Weyl groups, and as a straightforward
obvious analogy we emphasize the point set $L$ arising from the
application of $H^{\aff}_k$ groups in a similar way. Unlike the
crystallographic cases, the set $L$ covers densely the entire space
$\mathbb{R}^k$. Therefore some new elements have to be brought into 
consideration, which allow to select a suitable subset
$\Lambda\subset L$ for an application one  may have in mind. We point
out here three possible ways how that can be done, and pursue further
one of them.

The first one is the common cut-and-project method. It has been used
for many years and its construction does not require any visible
presence of affine Coxeter groups. Indeed, in this context $L$ arises as a result
of a projection of points of a higher dimensional crystallographic
lattice on a suitable subspace. That is particularly visible in the
algebraic definition of this projection (see \cite{CMP} and references
therein). The desired subset of $L$ is obtained by retaining only
the points which, under a complementary projection, fall into a bounded
`acceptance' region prescribed for $\Lambda$.  

The presence of an affine group becomes visible when only a finite number
of specifically affine transformations (translations) is required. An
example of such a case may be an algebraic formation of so called
carbon nanotubes and similar polytopes with non-spherical symmetry in
$\mathbb{R}^k$. Another example could be the modelling of concentric,
`onion-like', shell structures of carbon with $H_3$ symmetry
\cite{onion}. Obviously, every shell is one or several different
$H_3$-orbits, and transformations from shell to shell could
be provided by
$H^{\aff}_3$.

In this paper we pursue yet another way  in which $H^{\aff}_k$ can be exploited.  We
use it to get a finite subset $Q\subset \Lambda\subset L$, which lies in a
bounded region of the space $V$ and is a subset of a suitably defined fragment of a 
cut-and-project set. For that, we start from a seed point and allow no
more than a finite number $n<\infty$ of translations, while any number
of reflections from $H_k$ is admissible. Thus the value of
$n$ plays a similar role as does the acceptance window in the
cut-and-project case. In order to find their relation, we use the star
map \cite{CMP,P}, providing an explicit correspondence between the
points of the set $Q$ and the points of the corresponding 
`acceptance' window. 

The construction proposed here uses  the basic reflections of
$H^{\aff}_2$, $H^{\aff}_3$, and $H^{\aff}_4$ which are defined by
the simple roots encoded in the  extended Cartan matrices, or
equivalently in the extended Coxeter diagram. Any such extension
corresponds to a Coxeter group of infinite order. The latter are
obtained here similarly as in the framework of  Kac-Moody algebras
\cite{Kac:1985,Kass:1990}, where affine semisimple Lie algebras are
considered in parallel with the affine extensions of the corresponding
Weyl groups (crystallographic Coxeter groups). The extension allows
to identify a translation operation $T$ in $H^{\aff}_k$.
An iterative application of the basic reflections of $H_k$  and the operator 
$T\in H^{\aff}_k$ to a seed point in $\mathbb{R}^k$ then leads to  a
family of point sets $Q_k(n)$. The members of the families depend on an integer
valued {\sc cut-off parameter} $n$, which has simultaneously two
roles: (i) the value of $n$ determines the size of $Q_k(n)$, and
it plays a role similar to the acceptance window for
cut-and-project quasicrystals. In particular, it prevents
$Q_k(n)$ from becoming dense. We stress that for $n<\infty$, the
point sets
$Q_k(n)$ are of finite size, which distinguishes them from
cut-and-project quasicrystals which are generically infinite
structures.  A comparison with cut-and-project models shows
furthermore that
$H_k^{\aff}$-fragments $Q_k(n)$ are subsets of cut-and-project sets
with simply connected convex $H_k$-symmetric acceptance windows. 

The case of $H_2^{\aff}$-induced quasicrystals is
discussed in detail. Some properties, including bounds
on minimal next nearest neighbour distances
are discussed analytically. We remark, that combined
dilation-rotation symmetries for this type of point sets
have been investigated in \cite{Twarock:1999Wav} via
wavelet analysis. 

Since the technique for the derivation of the affine
extensions of noncrystallographic Coxeter groups is
similar to the one underlying the crystallographic 
case, it is instructive to briefly review the latter. We
will do this for the example of 
$sl(3)$, because like $H_2$, it has 2 simple roots of
the  same length. 

We then treat the noncrystallographic case indicating
explicitly the Cartan matrices for  $H^{\aff}_2$,
$H^{\aff}_3$ and $H^{\aff}_4$. We construct mathematical
models for $H^{\aff}_k$-induced fragments of quasicrystals with
particular emphasis on the case $k=2$, which is
investigated analytically and compared with the
cut-and-project scheme. 

\section{A short review of affine extensions based on the case of
$sl(3)$}
For convenience of the reader not familiar with the concept of
affine extensions, we briefly review the standard results for the
Weyl group of  $sl(3,\C )$ and of its affine extension. Our treatment
of the noncrystallographic case makes use of similar
considerations.  

The simple roots $\alpha_1$, $\alpha_2$ of  $sl(3,\C )$ span a real
Euclidean space $V$. All roots of $sl(3,\C )$ are of the same
length. In this section we adopt their normalization,
$(\alpha_k|\alpha_k)=2$, which is standard in Lie theory. Matrix
elements $a_{ij}$ of the Cartan matrix $A$ of
$sl(3,\C )$ are given in terms of their scalar products, 
\begin{equation}\label{Cartan matrix}
A:=(a_{ij})=
\left(\frac{2(\alpha_i\mid\alpha_j)}{(\alpha_j|\alpha_j)}\right)
= ((\alpha_i|\alpha_j)) 
=\left( \begin{array}{cc}
 2 & -1 \cr -1 & 2  
\end{array}\right)\,,\quad i,j=1,2.
\end{equation}
Simultaneously with
the basis of simple roots $\{\alpha_1,\alpha_2\}$, it is convenient to
work with the basis of fundamental weights $\{\omega_1,\omega_2\}$,
defined by
$$
(\alpha_j|\omega_k)=\tfrac12(\alpha_j|\alpha_j)\delta_{jk}
=\delta_{jk}\,.
$$
One has
\begin{equation}\label{bases}
\alpha_j=\sum_{k=1}^2a_{jk}\omega_k,\qquad
\omega_j=\sum_{k=1}^2(a^{-1})_{jk}\alpha_k\,,\qquad
A^{-1}=\tfrac13\left(\begin{array}{cc}
 2 & 1 \cr 1 & 2  
\end{array}\right)\,,
\end{equation}
where $(a^{-1})_{jk}$ are the matrix elements of the inverse
Cartan matrix $A^{-1}$.   Thus elements of the $j$th row of the 
Cartan matrix are the coordinates of the simple root $\alpha_j$  in
the $\omega$-basis, namely,
$\alpha_1=2\omega_1-\omega_2$ and $\alpha_2=-\omega_1+2\omega_2$. 

The  reflections $r_1$ and $r_2$, generating the Weyl group of
$sl(3,\C )$, act on a vector $v=v_1\omega_1+v_2\omega_2$ according to 
\begin{equation}\label{Weyl} 
r_j
v=v-\left(\frac{2(v|\alpha_j)}{(\alpha_j|\alpha_j)}\right)\alpha_j
=v-(v|\alpha_j)\alpha_j = v- v_j \alpha_j\,.
\end{equation} 
Due to the Weyl group symmetry of the root system, we can consider also 
reflections with respect to planes orthogonal to the other roots. 
 An affine extension $W_a$ of $W$ is obtained by introducing
the affine reflection $r_H^{\aff}$ as follows, where 
$\alpha_H=\alpha_1+\alpha_2=\omega_1+\omega_2$ is the highest root: 
\begin{equation}\label{affinereflect}
r_H^{\aff}v= v+\alpha_H-(v|\alpha_H)\alpha_H\,.
\end{equation}
From the particular cases $r_H^{\aff}0=\alpha_H$, 
$r_H^{\aff}\alpha_H=0$, $r_H^{\aff}\omega_1=\omega_1$, and 
$r_H^{\aff}\omega_2=\omega_2$, we see
that $r_H^{\aff}$ is a reflection in a plane orthogonal to
$\alpha_H$ and passing through the point $\frac12\alpha_H$ as well as
through the points $\omega_1,\ \omega_2$ rather than
through the origin. Consequently, 
$W_a$ contains the translation by $\alpha_H$, formed as the product of two 
reflections with respect to parallel mirrors:  
\begin{equation}
Tv := r_H^{\aff}r_H v 
= r_H^{\aff}\{ v - (v|\alpha_H)\alpha_H\} 
= v+\alpha_H\,.
\end{equation}

The extended Cartan matrix arises by adding to the simple roots
also the root $\alpha_0:=-\alpha_{H}$, using otherwise the same
conventions.  It leads to  
\begin{equation}
(a_{ij})=\left( \begin{array}{ccc}
 2 & -1 & -1 \\ -1 & 2 & -1 \\ -1 & -1 & 2
\end{array}\right)\,,\qquad i,j\in\{0,1,2\}\,.
\end{equation}
Note that such a matrix is subject to the following general
requirements: 
\begin{equation}\label{reqI}
a_{ii}=2\,,\quad
a_{ij}=a_{ji}\,,\quad
a_{ij}\in \integer^{\leq0}\, (i\not=j)\,,\quad
\det(a_{ij})=0 .
\end{equation}

Although the extended Cartan matrix cannot be inverted, one still could
define the dual basis independently in a similar way. In particular,
$$
\alpha_0 = 2\omega_0 - \omega_1 - \omega_2\,,\qquad
\alpha_1 = -\omega_0 +2\omega_1 - \omega_2\,,\qquad
\alpha_2 = -\omega_0 - \omega_1 +2\omega_2\,.
$$ 

The affine Weyl group operations $r_1$, $r_2$, and $T$ act on a
vector $v= v_1\omega_1 + v_2\omega_2$ according to 
\begin{equation}\label{ops}
\begin{aligned}
T v   & = v+\alpha_H &&= (v_1 + 1)\omega_1 + (v_2+1)\omega_2 
        &&= (v_1+1,v_2+1)\\
r_1 v & = v-(v|\alpha_1)\alpha_1 &&= -v_1\omega_1+(v_1+v_2)\omega_2)
        &&= (-v_1,v_1+v_2)\\
r_2 v & = v-(v|\alpha_2)\alpha_2 &&=(v_1+v_2)\omega_1-v_2\omega_2
        &&= (v_1+v_2,-v_2) \,.
\end{aligned}
\end{equation}
These transformations generate the  2-dimensional root lattice of
$sl(3,\C )$ from any single root or from zero. The reflections
$r_1$ and $r_2$ are subject to the defining identities of the Weyl group
of $sl(3, \C)$, see \eqref{2reflect} below.
The translation $T$ is not a cyclic operation, it can be repeated
any number of times. 

\section{The noncrystallographic Coxeter groups} 

 The Cartan matrices corresponding to the three
noncrystallographic Coxeter groups differ from the
crystallographic ones by the fact that their entries are
from the extension ring $\integer[\t]:=\lbrace
a+\t b\mid a,b\in\integer\rbrace$, where the irrationality
 is the golden mean 
$$
\t :=\half(1+\sqrt{5})\,,\qquad
\t':=\half(1-\sqrt{5})=1-\t=-\frac1{\t}\,.
$$
 Thus the conditions
(\ref{reqI}) on the extended Cartan matrices become  
\begin{equation}\label{reqII}
a_{ii}=2\,,\quad
a_{ij}=a_{ji}\,,\quad
a_{ij}\in \integer[\t]^-:=\lbrace
     x\in\integer[\t]\mid x\leq0\rbrace\,,\quad
\det(a_{ij})=0 .
\end{equation}

Introducing again the additional root via $\alpha_0=-\alpha_H$ where
$\alpha_H$ is the highest root, the extended Cartan matrices are 
obtained from the Cartan matrices of $H_j$, $j=2,3,4$. The
corresponding groups will be denoted as $H^{\aff}_j$, $j=2,3,4$, 
respectively. A direct calculation shows that such matrices
are the unique ones fulfilling all requirements (\ref{reqII}). We stress
that the condition $a_{ij}\in \integer[\t]^-$ is crucial for
uniqueness.  Without it, several Cartan matrices can be found which
fulfill all other conditions. Such matrices are shown in Appendix A.

We now discuss the three cases separately. We
describe the $H_2$ case in details. The other two, $H_3$ and $H_4$, are exact
analogies. There we provide important steps and the result of the
considerations.

Furthermore, note that it is not possible to obtain $H^{\aff}_j$,
$j=2,3,4$ via a projection from a group with $2(k+1)$ simple roots,
which has a Cartan matrix obeying (\ref{reqI}); it would be necessary
to relax the third assumption and admit also positive entries
$a_{ij}$. 

Unlike the crystallographic case, we normalize the simple roots of $H_k$
to be of length one. Note that Cartan matrices do not depent on root
normalization.
\bigskip

\subsection{The case of $H^{\aff}_2$ as an affine extension of $H_2$.}
\par
The Coxeter group $H_2$ is isomorphic to the dihedral group of order
10 and its  root system  can be modeled in the complex plane by the 
10th roots of unity. The root system $\Delta_2$ is the  union of the sets
of positive and negative roots. Choosing the simple roots as
$\alpha_1=1$ and $\alpha_2=\exp (\frac{4\pi i}{5})$, the roots
\begin{equation}\label{rootsH2}
\Delta_2 =\{ 
\pm\alpha_1,\ \pm\alpha_2,\ \pm(\alpha_1+\t \alpha_2),\
\pm(\t \alpha_1+ \alpha_2),\ \pm(\t\alpha_1+\t
\alpha_2)\}
\end{equation}
form the vertex set of a regular decagon inscribed into the unit circle.
Now $(\alpha|\alpha)=1$ for any $\alpha\in\Delta_2$.
From (\ref{Cartan matrix}) we find the Cartan matrix and its
inverse,  
\begin{equation}\label{2by2CM}
A=\left(\frac{2(\alpha_i|\alpha_j)}{(\alpha_j|\alpha_j)}
\right)= 2((\alpha_i|\alpha_j))=\left(
\begin{array}{cc}
 2 & -\t \cr -\t & 2  
\end{array}\right)\,,\quad 
A^{-1}=\frac1{3-\t}\left(\begin{matrix} 2&\t\\\t&2\end{matrix}\right).
\end{equation}
Here, as before, $\omega_1$ and $\omega_2$ are the basis vectors of the
$\omega$-basis defined by $2(\alpha_j|\omega_k)=\delta_{jk}$. It
follows that
$$\begin{aligned}
\alpha_1  &= 2\omega_1-\t\omega_2\,,\qquad 
&\alpha_2 &= -\t\omega_1+2\omega_2\,;\\
\omega_1  &= \tfrac1{3-\t}(2\alpha_1+\t\alpha_2)\,,\qquad
&\omega_2 &= \tfrac1{3-\t}(\t\alpha_1+2\alpha_2)\,.
\end{aligned} 
$$
The highest root is  $\alpha_{H}=\t(\alpha_1+\alpha_2)
= -\t'(\omega_1+\omega_2)$.

Taking the extension root as $\alpha_0:=-\alpha_{H}$ and letting the
indices in (\ref{2by2CM}) take the values $0$, $1$, and $2$, the Cartan
matrix of the affine extension $H^{\aff}_2$ and its simple roots in the 
$\omega$-basis turn out to be
\begin{equation}\label{GenCart}
\begin{pmatrix}
 2&\t'&\t'\\ \t'&2&-\t\\ \t'&-\t&2
\end{pmatrix}\,;\qquad\qquad
\begin{aligned}
\alpha_0 &= 2\omega_0+\t'\omega_1+\t'\omega_2\,,\\
\alpha_1 &= \t'\omega_0+2\omega_1-\t\omega_2\,,\\
\alpha_2 &= \t'\omega_0-\t\omega_1+2\omega_2\,.
\end{aligned}
\end{equation}
Indeed, using $2(\alpha_0|\alpha_1) =2(-\alpha_H|\alpha_1) =
2(-(\t\alpha_1+\t \alpha_2)|\alpha_1) 
=-2\t +\t^2= \t'$  in
(\ref{2by2CM}), we get the matrix elements of (\ref{GenCart}). 

\bigskip
The corresponding Coxeter diagram is given in Fig. 1. 

\begin{figure}[ht]\label{one}
\begin{center}
\includegraphics[width=1.8cm,keepaspectratio]{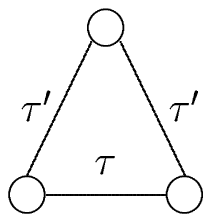}\\*
\end{center}
\caption{$H^{\aff}_2$--diagram}
\end{figure}

\bigskip
The nodes of the diagram stand for the simple roots. A direct link between two
nodes indicates that  the two roots are not orthogonal in the Euclidean plane
spanned by them.  Two roots are orthogonal if there is no direct link between
the corresponding nodes. The label attached to a link is determined by the
off-diagonal matrix element of the Cartan matrix: no label is shown if such
element is $-1$; if it is $-\t$, the link is labeled by $\t$;
if the matrix element is $\t'$, the label is
$\t'$. 

The reflection $r_0$ with respect to the plane orthogonal to
$\alpha_H$ and passing through the origin is defined by the general
formula
\eqref{Weyl}, where now one has to use the roots of $H_k$
rather than those of $sl(3,\C )$. Similarly, the affine reflection
$r_H^{\aff}$  is defined as the reflection in the plane orthogonal to
$\alpha_H$ and passing through the point $\tfrac12\alpha_H$. Due to
the different normalization of the bases $\{\alpha_j\}$ and
$\{\omega_k\}$, some modification appears in the corresponding
formulas. Thus instead of (\ref{affinereflect}), we now have
$$
r_H^{\aff}v = v + \left\{1 - 2(v|\alpha_H)\right\}\alpha_H\,.
$$
Note the instructive particular cases  
$$
r_H^{\aff}0 = \alpha_H, \quad 
r_H^{\aff}\alpha_H = 0, \quad\text{and}\quad
r_H^{\aff}\left(\tfrac1\t\omega_j\right)
                   =\left(\tfrac1\t\omega_j\right), \quad j=1,2\,.
$$
The product of the reflections $r_H^{\aff}r_0$ is the translation
operator $T$.   
\begin{equation}
\begin{aligned}
Tv = r_H^{\aff}r_0 
   &= r_H^{\aff}\left\{(v-2(v|\alpha_H)\alpha_H)\right\} \\
   &= v-2(v|\alpha_H)\alpha_H + \left\{1 -
       2(v-2(v|\alpha_H)\alpha_H|\alpha_H)\right\}\alpha_H \\
   &= v+\alpha_H\,.
\end{aligned}
\end{equation}
Explicitly in the $\omega$-basis, we have 
\begin{equation}\label{trafos2}
\begin{aligned}
T v   & = v+\alpha_H 
        &&= (v_1 - \t')\omega_1 + (v_2-\t')\omega_2
        &&=(v_1-\t',v_2-\t')\\
r_1 v & = v-2(v|\alpha_1)\alpha_1 
        &&= -v_1\omega_1+(\t v_1+v_2)\omega_2
        &&=(-v_1,\t v_1+v_2)\\ 
r_2 v & = v-2(v|\alpha_2)\alpha_2 
        &&=(v_1+\t v_2)\omega_1-v_2\omega_2 
        &&=(v_1+\t v_2,-v_2)\,.
\end{aligned}
\end{equation}

The reflections are subject to the $H_2$ group identities (\ref{2reflect}), namely,  
\begin{equation}
r_1^2=r_2^2=1,\qquad (r_1r_2)^5=1\,.
\end{equation}
In contrast to $r_1$ and $r_2$, the translation $T$ can be repeated any number of times without producing the same points.

Below, these transformations are used in order to build 2-dimensional quasicrystalline
point sets similarly as in the previous case for the $sl(3,\C )$ root lattice. 
Straightforwardly repeated applications of the three operations
in every possible sequence, without further restrictions, would produce
a dense point set covering the whole plane. Note that if the coordinates
$v_1,\ v_2$ of the seed point
$v$ are in
$Z[\t]$, every point of the set has its coordinates in
$Z[\t]$.
\bigskip

\subsection{The case of $H^{\aff}_3$ as an affine extension of $H_3$}
\par
The root system of $H_3$ consists of 30 roots; they can be found in
\cite{CKPS}. Relative to an orthonormal basis, normalized to length $1$
rather than $\sqrt2$ like in the preceeding subsection, they can be
modeled as 
\begin{equation}\label{icosH3}
\Delta_3 = \left\lbrace{
\begin{array}{cl}
(\pm 1,0,0) & \mbox{ and all permutations }\\
\half(\pm 1,\pm \t',\pm \t) & \mbox{ and all even permutations }
\end{array}
}\right\rbrace\,.
\end{equation}
Geometrically, the root polytope of  $H_3$ is formed by 12 equilateral pentagons and 20 equilateral triangles. It has 30 vertices given by the elements in $\Delta_3$ and 60 edges.  
It is possible and sometimes advantageous to consider the roots of
$\Delta_3$ given in \eqref{icosH3} as purely imaginary
quaternions of special kind, called icosians \cite{Moody:1993,CMP}.

A possible choice of simple roots in the orthonormal basis is
$$
\alpha_1=(0,0,1)\,,\quad\alpha_2=\tfrac12(-\t',-\t,-1)\,,\quad
\alpha_3=(0,1,0)\,.
$$
The Cartan matrix of $H_3$ and its inverse,
$$
A=
\begin{pmatrix}
 2 &    -1 & 0 \cr 
-1 &    2  & -\t \cr
 0 & -\t &   2  
\end{pmatrix}\,,\qquad
A^{-1}=\frac12
\begin{pmatrix}
 2+\t  & 2+2\t & 1+2\t \cr 
 2+2\t & 4+4\t & 2+4\t \cr
 1+2\t & 2+4\t & 3+3\t  
\end{pmatrix}\,,
$$
are used to find
\begin{alignat*}{2}
\alpha_1 &= 2\omega_1-\omega_2\,,  &\qquad 
\omega_1 &=\tfrac12((2+\t)\alpha_1+2\t^2\alpha_2+\t^3\alpha_3)\,,\\
\alpha_2 &= -\omega_1+2\omega_2-\t\omega_3\,, &\qquad 
\omega_2 &= \t^2\alpha_1+2\t^2\alpha_2+\t^3\alpha_3\,,\\
\alpha_3 &= -\t\omega_2+2\omega_3\,, &\qquad 
\omega_3 &=\tfrac12(\t^3\alpha_1+2\t^3\alpha_2+3\t^2\alpha_3)\,.
\end{alignat*}
The highest root is $\alpha_H=
\t\alpha_1+2\t\alpha_2+\t^2\alpha_3= -\t'\omega_2 = (1,0,0)$.

The affine extension of the Cartan matrix, $H^{\aff}_3$, and the
simple roots in the $\omega$-basis are
$$
\begin{pmatrix}
 2    &  0  & \t' & 0     \\ 
 0    &  2  &  -1   & 0     \\
\t' &  -1 &   2   & -\t \\
 0    &  0  & -\t &   2  
\end{pmatrix}\,; \qquad
\begin{aligned}
\alpha_0 &= 2\omega_0+\t'\omega_2\\
\alpha_1 &= 2\omega_1- \omega_2\\
\alpha_2 &= \t'\omega_0- \omega_1+2\omega_2-\t\omega_3\\
\alpha_3 &= -\t\omega_2+2\omega_3\,.
\end{aligned}
$$
This corresponds to the graph in Fig. 2. 
\medskip

\begin{figure}[ht]\label{two} 
\begin{center}
\includegraphics[width=2.8cm,keepaspectratio]{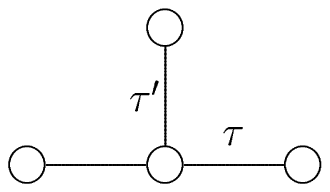}\\*
\end{center}
\caption{$H^{\aff}_3$--diagram}
\end{figure}
\noindent

The reflections $r_1$, $r_2$, $r_3$ as well as $T$ are built from
the general expressions \eqref{Weyl} and \eqref{T}, where the
indices run through three values and the roots of $H_3$ are used.
They act on
$v=(v_1,v_2,v_3)$ in the 
$\omega$-basis according to 
\begin{equation}\label{trafos3}
\begin{aligned}
T v   & = v+\alpha_H &&=(v_1,v_2-\t', v_3)\\
r_1 v & = v-2(v|\alpha_1)\alpha_1  &&=(-v_1, v_1+v_2,v_3)\\
r_2 v & = v-2(v|\alpha_2)\alpha_2  &&= (v_1+v_2, -v_2, v_3+\t v_2)\\ 
r_3 v & = v-2(v|\alpha_3)\alpha_3  &&= (v_1, v_2+\t v_3,-v_3)\,.
\end{aligned}
\end{equation}
\bigskip

\subsection{The case of $H^{\aff}_4$ as an affine extension of $H_4$}

The root system $\Delta_4$ of $H_4$ contains 120 roots, they are
found in \cite{CKPS} in terms of simple roots. They can be modeled
\cite{CMP} as the set, 
\begin{equation}
\Delta_4 = \left\lbrace{
\begin{array}{cl}
\half(\pm 1,\pm 1,\pm 1,\pm 1), (\pm 1,0,0,0) & \mbox{ and all
permutations }\\
\half(0,\pm 1,\pm \t',\pm \t) & \mbox{ and all even permutations }
\end{array}
}\right\rbrace\,
\end{equation}
in an orthonormal basis, or equivalently as quaternions \cite{Moody:1993,CMP}. 
Equipped with quaternionic multiplication,
they stand for the elements of the icosahedral group.

As simple roots, one may choose 
\begin{equation}
\begin{aligned}
\alpha_1   &=\tfrac12(-\t',-\t,0,-1)\,,\qquad
 && \alpha_2=\tfrac12(0,-\t',-\t,1)\,,\\
\alpha_3   &=\tfrac12(0,1,-\t',-\t)\,, 
 && \alpha_4=\tfrac12(0,-1,-\t',-\t)\,.
\end{aligned}
\end{equation}
The highest root of $H_4$ is then
$\alpha_H= 
2\t\alpha_1+ \sqrt5\t^2\alpha_2+ 2\t^3\alpha_3+ \t^4\alpha_4
= -\t'\omega_1=(1,0,0,0)$. The $H_4$-Cartan matrix and its
inverse are as follows: 
$$
A=
\left( \begin{matrix}
 2 & -1 &  0  &  0  \\ 
-1 &  2 & -1  &  0  \\
 0 & -1 &  2  & -\t \\
 0 &  0 & -\t &   2  
\end{matrix}\right)\,,\qquad
A^{-1}=
\left( \begin{matrix}
2+2\t & 3+4\t  & 4+6\t   & 3+5\t     \\ 
3+4\t & 6+8\t  & 8+12\t  & 6+10\t     \\
4+6\t & 8+12\t & 12+18\t & 9+15\t \\ 
3+5\t & 6+10\t & 9+15\t  & 8+12\t  
\end{matrix}\right)\,.
$$
As generalized Cartan matrix and simple roots in the $\omega$-basis, we
obtain 
$$
\left(\begin{matrix}
 2    &  \t' &  0  &   0   &  0      \\ 
\t' &    2   &  -1 &   0   &  0      \\
 0    &   -1   &   2 &   -1  &  0      \\
 0    &    0   &  -1 &    2  & -\t   \\ 
 0    &    0   &  0  & -\t &  2
\end{matrix}\right)\,,\qquad
\begin{aligned}
\alpha_0 &= \tfrac12(2\omega_0+\t'\omega_1)\,,\quad
&&\alpha_1 = \tfrac12(\t'\omega_0+2\omega_1- \omega_2)\,,\\
\alpha_2 &= \tfrac12(-\omega_1+2\omega_2-\omega_3)\,,\qquad
&&\alpha_3 = \tfrac12(-\omega_2+2\omega_3-\t\omega_4)\,,\\
\alpha_4 &= \tfrac12(-\t\omega_3+2\omega_4)\,,
\end{aligned}
$$
and the corresponding Coxeter diagram is depicted in Fig. 3. 

\begin{figure}[ht]\label{three}
\begin{center}
\includegraphics[width=4.5cm,keepaspectratio]{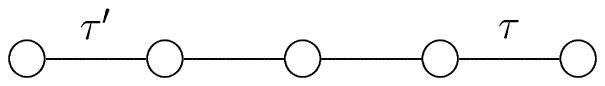}\\*
\end{center}
\caption{$H^{\aff}_4$--diagram}
\end{figure}

The translation and the four reflections act in 4-space on a
point $v=(v_0,v_1,v_2,v_3,v_4)$ with coordinates in the $\omega$-basis according to: 
\begin{equation}\label{trafos4}
\begin{aligned}
T v   & = (v_1-\t', v_2, v_3, v_4)\\
r_1 v & = (-v_1, v_1+v_2, v_3, v_4)\\
r_2 v & = (v_1+v_2, -v_2, v_2+v_3, v_4)\\ 
r_3 v & = (v_1,v_2+v_3, -v_3, v_4+\t v_3)\\ 
r_4 v & = (v_1, v_2, v_3+\t v_4, -v_4)\,.
\end{aligned}
\end{equation}
Here, as in \eqref{ops},\eqref{trafos2},\eqref{trafos3}, the
reflections are cyclic operations of order two, while the
translation $T$ can be repeated any number of times. Products of
two reflections are rotations around the origin. Their order is
determined by the matrix elements of the corresponding
Cartan matrix
\begin{equation}\label{2reflect}
(r_jr_k)^M= 1 \quad\text{where}\quad
\left\{\begin{matrix}   
M=1 \quad &\text{if}\ &a_{jk}&= 2\\
M=2 \quad &\text{if}\ &a_{jk}&= 0\\
M=3 \quad &\text{if}\ &a_{jk}&=-1\\
M=5 \quad &\text{if}\ &a_{jk}&=-\t,\t'
\end{matrix}\right.
\end{equation}
\bigskip

\section{Construction of $H^{\aff}$-induced quasicrystals} 

In this section, we illustrate an application of the transformations
(\ref{trafos2}), (\ref{trafos3}) and (\ref{trafos4}) in order to
generate point sets in Euclidean spaces of dimensions 2, 3, and 4, 
respectively, which resemble fragments of quasicrystals. More precisely, the number 
of allowed translations plays a similar role as the acceptance window of a cut-and-project quasicrystal and 
a certain neighbourhood of the seed point contains all the points of such a quasicrystal and only at the periphery the fragment has fewer points. 
The idea ofthe construction is to use the reflections on a given seed point in
every possible way, 
while using
the translation $T$ only for a fixed finite number of times. 

We start by describing the construction in detail for the case of
$H_2^{\aff}$, after that the other two cases are straightforward. 

Transformations (\ref{trafos2}) acting in 2-space can be represented
using $2\times2$ matrices. For notational convenience, we use the symbol
$v=v_1\omega_1+v_2\omega_2$ also for the column matrix $(v_1\
v_2)^T$: 
\begin{equation*}\label{T}
Tv=
 \begin{pmatrix} 1 & 0 \\ 0 & 1 \end{pmatrix} v - \t' 
\begin{pmatrix} 1\\ 1 \end{pmatrix}
= (v_1-\t')\omega_1 + (v_2-\t')\omega_2\,,
\end{equation*}
\begin{equation}\label{R1}
r_1 v:=\begin{pmatrix} -1 & 0 \\ \t & 1  \end{pmatrix} v
= -v_1\omega_1 + (\t v_1 + v_2)\omega_2\,,
\end{equation}
\begin{equation*}\label{R2}
r_2 v:=\begin{pmatrix} 1 & \t \\ 0 & -1  \end{pmatrix} v
= (v_1+\t v_2) - v_2\omega_2\,.
\end{equation*}
A straightforward calculation shows that also
$$
(r_1r_2)^5=
\left(
\begin{pmatrix}-1&0 \\ \t& 1 \end{pmatrix}
\begin{pmatrix} 1&\t\\ 0 &-1 \end{pmatrix}
\right)^5
=\begin{pmatrix}-1&-\t\\ \t&\t\end{pmatrix}^5
=\begin{pmatrix} 1&0\\0 &1\end{pmatrix}\,.
$$

In order to change to
Cartesian coordinates, a further transformation is needed: 
\begin{equation}\label{S}
x:=\left( \begin{array}{cc}
 0 & r \\ 1 & \frac{\t}{2}   
\end{array}\right) v\,,\qquad r=\sqrt{1-\frac{\t^2}{4}}\,.
\end{equation}

Note that the transformations $r_1$ and $r_2$ act as reflections at the mirrors
perpendicular to the simple roots of $H_2$, which are collinear with the
$\omega$-basis and intersecting at the origin.
Their relative angle is $2\pi/10$ and $T$ defines a translation along
their bisector.

The iterate action of the transformations $T$, $r_1$ and $r_2$ in
arbitrary order starting from the origin leads to a point set which
fills the plane densely after an infinite number of iterations. If the
iteration is stopped after a finite number of steps,  a discrete point
set is obtained. The application of $r_1$ and $r_2$ transforms a point
within the same  $H_2$-orbit. They are equidistant from the origin.
Points translated by $T$ are on different
$H_2$-orbits. All the points generated by the three operators from one
seed point are within one $H^{\aff}_2$-orbit.

\begin{de}
A point $v$ is said to be dominant precisely if its coordinates in the
$\omega$-basis are non-negative.
\end{de}

It is convenient to characterize an orbit of $H_k$ by its unique point
({\it dominant point}) which is the only one in its orbit  with
non-negative coordinates in the $\omega$-basis. It is easily recognizable by
this property. Since dominant points encode the information about the whole 
$H_k$-orbit, they are a useful tool for the construction and analysis of the 
point sets. 

The size of an $H_2$-orbit is readily found from its dominant representative
by the following rule:
\begin{equation}
\begin{array}{llllllll}
&{\rm orbit\ size}\quad &10&:\quad && v=(a\ b)^T\,,\quad && a>0,b>0 
\notag\\
&{\rm orbit\ size}\quad & 5&:\quad && v=(a\ 0)^T\quad{\rm or}\quad v=(0\
b)^T\,,\quad && a>0,b>0 
\\
&{\rm orbit\ size}\quad & 1&:\quad && v=(0\ 0)^T\,.\quad && 
\notag
\end{array}
\end{equation}
Here, $(a b)^T$ denotes the transposition of the row matrix $(a b)$. 

For example, applying $r_1$ and $r_2$ to $v=(a\ 0)^T$, $a>0$, according
to (22), one gets the following five points of the $H_2$ orbit, which correspond to 
the vertices of an equilateral pentagon centered at the origin: 
\begin{equation}\label{penta}
\begin{pmatrix} a\\0 \end{pmatrix},\quad
\begin{pmatrix} -a\\a\t \end{pmatrix},\quad
\begin{pmatrix} a\t\\-a\t \end{pmatrix},\quad
\begin{pmatrix} -a\t\\a \end{pmatrix},\quad
\begin{pmatrix} 0\\-a \end{pmatrix}\,.
\end{equation}
Similarly starting from the point $v=(0\ a)^T$, $a>0$, one gets
another $5$-point orbit consisting of the negatives of \eqref{penta}.
Any further application of
$r_1$ and
$r_2$ would bring no new point. 

The square length of the vectors given in the $\omega$-basis is
calculated using the  inverse of the Cartan matrix (\ref{2by2CM}):
$$
2(a\omega_1+b\omega_2|a\omega_1+b\omega_2)=
\frac1{3-\t}\begin{pmatrix} a &b \end{pmatrix}\begin{pmatrix} 2 & \t\\
\t & 2 
\end{pmatrix}
\begin{pmatrix} a\\b \end{pmatrix} = \frac{2(a^2 + ab\t +
b^2)}{3-\t}\,.
$$

\begin{de}\label{affQCdef} 
Let $O$ denote the origin of coordinates, and let $s^m(T,r_1,r_2)$
denote the set of all words formed by the letters $T$, $r_1$ and $r_2$
in which $T$ appears precisely $m$ times. The set of points 
\begin{equation}\label{affQC} 
Q_2(n):=\{s^m(T,r_1,r_2)O\mid m\leq n\}
\end{equation} 
is called an $H_2^{\aff}$-induced quasicrystal fragment; $n$ is the
cut-off-parameter. 
\end{de}
\noindent
Due to the identities (\ref{2reflect}), the point set $Q_2(n)$ is
finite and $H_2$-symmetric with respect to the origin. More precisely, the fact that we allow an arbitrary number of actions of $r_1$ and $r_2$ after each translation enforces the finite patches to have circular boundaries. We make this assumption here in view of applications, because this is the situation one encounters e.g. for carbon onions in the study of fullerenes, or, this is what we expect for the growth process of a quasicrystal fragment which is not exposed to any particular obstacles. We remark that it is possible to change (\ref{affQC}) by requiring that after the last translation, no further actions of $r_1$ or $r_2$ take place. In this case, $H_2$-symmetry with respect to the origin would no longer be present in the model.  

Note that due to $H_2$-symmetry, each $Q_2(n)$ can be decomposed into concentric shells containing all the points at
the same distance from the origin. In general each shell is a union of
several decagons and pentagons, except for the origin which alone is a
one-point shell. The outer most shell of 
$Q_2(n)$ is  the equilateral decagon with dominant point $T\dots TO$, where 
$n$ translation operators $T$ are applied to the origin.
\medskip

Clearly $Q_2(0)$ is just one point, the dominant point $O=(0,0)$. 
The set $Q_2(1)$ contains the origin and the vertices of the
decagon of $H_2$ roots. Among the latter the highest root $TO$ is the dominant
one. Thus $Q_2(1)$ contains eleven points. It is the union of $Q_2(0)$
and the orbit of the roots of $H_2$:
\begin{equation}
\begin{array}{ll}
&O,\qquad TO=\t\alpha_1+\t\alpha_2,\qquad 
r_1TO=\t\alpha_1+\alpha_2,\qquad 
r_2TO=\alpha_1+\t\alpha_2,
\\ 
&r_1r_2TO=\alpha_1,\qquad 
r_2r_1TO=\alpha_2,\qquad 
r_1r_2r_1TO=-\alpha_2,\qquad 
r_2r_1r_2TO=-\alpha_1,
\\ 
&r_2r_1r_2r_1TO=-\t\alpha_1-\alpha_2,\qquad 
r_1r_2r_1r_2TO=-\alpha_1-\t\alpha_2,
\\ 
&r_1r_2r_1r_2r_1TO=r_2r_1r_2r_1r_2TO=-\t\alpha_1-\t\alpha_2\,.
\end{array}
\end{equation}
The equality of the words here and the absence of words involving
$r_j^2$ are consequences of the defining identities (\ref{2reflect}) of
the group $H_2$. Further applications of $r_1$ and $r_2$  yield no new
points. 

The set $Q_2(2)$ is obtained by shifting $Q_2(1)$ by $T$, i.e. by
the highest root  $\alpha_H$, and by applying to the result all possible
$r$'s. It contains 
61 distinct points. 

It decomposes into the sum of four  orbits of 10 points with the
dominant points 
\begin{equation}
\begin{array}{rclrcl}
2 \alpha_H & = & 2 \t (\alpha_1 + \alpha_2) & 
\t \alpha_H  & = &  \t^2 (\alpha_1 + \alpha_2) \\
\alpha_H & = & \t (\alpha_1 + \alpha_2) &
(-\t') \alpha_H & = & \alpha_1 + \alpha_2\,,
\end{array}
\end{equation}
 four orbits of 5 points with dominant points
\begin{equation}
\begin{array}{rclrcl}
r_1 (-\alpha_1+\alpha_H) & = & 2 \alpha_1 + \t \alpha_2 & 
r_2 (-\alpha_2 + \alpha_H)  & = &  2 \alpha_2 +\t \alpha_1 \\
\t^2 \alpha_1 + 2 \t \alpha_2 & &  &
\t^2 \alpha_2 + 2 \t \alpha_1 \,,& & 
\end{array}
\end{equation}
and the origin.

Compare with Fig. 4, where $Q_2(2)$ is depicted. 

\begin{figure}[ht]\label{four}
\begin{center}
\includegraphics[width=5.0cm,keepaspectratio]{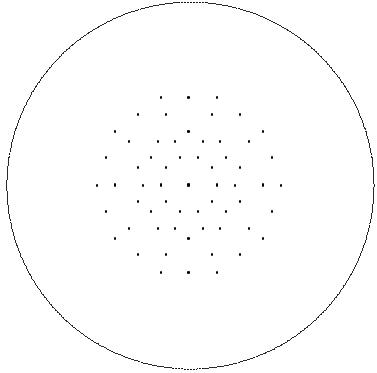}\\*
\end{center}
\caption{The point set $Q_2(2)$}
\end{figure}

\medskip

In our construction $n$ plays a similar role as the acceptance windows for
cut-and-project quasicrystals, because it ensures that instead of a dense
set, a discrete point set is obtained. It is therefore interesting to
identify $Q_2(n)$ as a subset of points of a cut-and-project set as far as
possible. 

Recall that a cut-and-project point set is completely determined by its
acceptance window. There is a 1--1 map between the points of the window
and of the cut-and-project set. A finite fragment of an (infinite)
cut-and-project point set allows to determine its acceptance window
only within certain bounds \cite{Patera:1998c}. The larger is the fragment, the tighter
are the bounds on the window. All the points of 
$Q_2(n)$ are found inside of the decagon formed by the outermost shell
whose points are at the distance $n\t|\alpha_1+\alpha_2|$ from the
origin, where $\t|\alpha_1+\alpha_2|$ corresponds to the length of
any $H_2$ root.

Note that $Q_2(n)$ is invariant under 10-fold rotational symmetry by
construction. It would also be possible to define
an aperiodic point set based on the operations $T$, $r_1$ and $r_2$ which
does not have this property. One may e.g. take instead of $s^r(T,r_1,r_2)$
all sequences which end after the operation $T$ in order to break this
symmetry. 

\section{Generalization to $H_3^{\aff}$ and $H_4^{\aff}$} 

An extension of the previous construction to 3 and 4 dimensions is
straightforward. In analogy to the previous section, we obtain the
following  operators $T$ and $r_j$ from (\ref{trafos3}) and 
(\ref{trafos4}). 

\begin{enumerate}

\item For $H_3^{aff}$ and a vector $v=(v_1\ v_2\ v_3)^T$ with coordinates in the $\omega$-basis we have:

\begin{equation}\label{Tb}
T v:=\left( \begin{array}{ccc}
 1 & 0 & 0\\ 0 & 1 & 0\\ 0 & 0 & 1 
\end{array}\right)v -\t' 
\left( \begin{array}{c}
 0\\ 1\\ 0 
\end{array}\right)
=v_1 \omega_1 + (v_2-\t') \omega_2 + v_3 \omega_3
\,,
\end{equation}
\begin{equation}\label{R1b}
r_1 v=\begin{pmatrix}
 -1 & 0 & 0\\ 1& 1 & 0\\ 0 & 0 & 1
\end{pmatrix} v\,,
\quad
r_2 v=\begin{pmatrix}
 1 & 1 & 0 \\ 0 & -1 & 0 \\ 0 & \t & 1 
\end{pmatrix} v\,,
\quad
r_3 v=\begin{pmatrix}
 1 & 0 & 0 \\ 0 & 1 & \t \\ 0 & 0 &-1
\end{pmatrix} v\,.
\end{equation}

\item For $H_4^{aff}$ and a vector $v=(v_1\ v_2\ v_3\ v_4)^T$ with coordinates in the $\omega$-basis we have:

\begin{equation}\label{Tc}
T(v_0) v:=\left( \begin{matrix}
 1 & 0 & 0 & 0\\ 0 & 1 & 0 & 0 \\ 0 & 0 & 1 & 0\\ 0 & 0 & 0 & 1
\end{matrix}\right)v - \t' 
\left( \begin{matrix}
 1\\ 0 \\ 0\\ 0 
\end{matrix}\right)
= (v_1-\t')\omega_1 + v_2 \omega_2 +  v_3 \omega_3
+ v_4 \omega_4
\,,
\end{equation}
\begin{equation}\label{R1c}
r_1 v:=\left(\begin{matrix}
 -1 & 0 & 0 & 0\\ 1& 1 & 0 & 0\\ 0 & 0 & 1 & 0\\ 0 & 0 & 0 & 1
\end{matrix}\right) v\,,
\quad 
r_2 v:=\left( \begin{matrix}
 1 & 1 & 0 & 0\\ 0 & -1 & 0 & 0 \\ 0 & 1 & 1 & 0\\ 0 & 0 & 0 & 1
\end{matrix}\right) v\,,\quad
\end{equation}
\begin{equation}\label{R3c}
r_3 v:=\left( \begin{matrix}
 1 & 0 & 0 & 0\\ 0 & 1 & 1 & 0 \\ 0 & 0 &-1 & 0\\ 0 & 0 & \t & 1
\end{matrix}\right) v\,,
\quad 
r_4 v:=\left(\begin{matrix}
 1 & 0 & 0 & 0\\ 0 & 1 & 0 & 0 \\ 0 & 0 & 1 & \t\\ 0 & 0 & 0 & -1
\end{matrix}\right) v\,,
\end{equation}

\end{enumerate}

In both cases, the orbit sizes can be found in \cite{CKPS}.

\begin{de}\label{affQdef}
Let $s^m(T,r_1,...,r_k)$ denote the set of all sequences formed by the
operators $T$ and $r_1,\ldots,r_k$ in which $T$ appears precisely $m$
times; O denotes the origin of coordinates. Then 
\begin{equation}\label{affQC34}
 Q_k(n):=\lbrace s^m(T,r_1,...,r_k)O \mid
m\leq n\rbrace
\end{equation} 
is called $H_k^{\aff}$-induced quasicrystal fragment for
$k=3,4$; $n$ is the cut-off-parameter. 
\end{de} Note that $Q_k(n)$ describes a $k$-dimensional point set. 

\section{Investigation of $H_2^{\aff}$-induced quasicrystals $Q_2(n)$}

In this section we analyze the point sets $Q_2(n)$ which have been defined in Def. \ref{affQCdef}. Recall that the set $Q_2(n)$ was obtained 
from the origin via an application of three operations, the translation $T$ and 
the reflections $r_1$ and $r_2$ subject to the condition that the operation $T$ 
occurs precisely $k$ times, whereas any number of reflections is permitted. 
As pointed out before, the point set is characterized by its dominant points. 

\subsection{The dominant points $(a b)^T$ with $a=b$:}

We start by an investigation of the dominant points $(a\ b)^T$ with $a=b$. Note that they are given by multiples of the highest root $\alpha_H=\tau(\alpha_1+\alpha_2)$ and are thus located on the line $\R \alpha_H$. We introduce the notation $L_{\alpha_H}(n)$ for the finite point set  given as the intersection of the 2-dimensional point set $Q_2(n)$ with the line $\R \alpha_H$:
\begin{equation}\label{Ldef}
L_{\alpha_H}(n) := \R \alpha_H\cap Q_2(n)\,.
\end{equation}
The first step in our analysis of $Q_2(n)$ will be a description and analysis of $L_{\alpha_H}(n)$. 

For the remainder of this paper, we model the root system of $H_2$ in (\ref{rootsH2})
in the complex plane by $\xi^j$ where 
\begin{equation}\label{rootpoints}
\xi := \exp i\frac{\pi}{5}   \,,
\end{equation}
that is $\xi^0,\ldots,\xi^9$ number the roots of unity anticlockwise starting
from $\xi^0=1$. 
We remark that due to the 10-fold rotational symmetry of $Q_2(n)$, the set $L_{\alpha_H}(n)$ in (\ref{Ldef}) coincides -- when viewed as a 1 dimensional point set -- with the sets 
\begin{equation}
L_{\xi^j}(n) := \R \xi^j\cap Q_2(n) \mbox{ for } \xi^j\in\Delta_2\,.
\end{equation}
Thus, the results obtained for any $L_{\xi^j}(n)$ with $\xi^j\in\Delta_2$ translate immediately into each other. 

\subsubsection{Description of $L_{\alpha_H}(n)$:}

\noindent
We start by expressing the points in $Q_2(n)$ in a more convenient way. 
For this, recall that by definition  
\begin{equation}\label{points}
x\in Q_2(n) \Leftrightarrow x=
{R}_l T {R}_{l-1} T \ldots T {R}_1 T  O, \qquad l\leq n
\end{equation}
where $O$ denotes the origin and ${R}_j$ for $j=1,\ldots,n$ denotes 
a product of basic reflections $r_1$ and $r_2$, i.e. an element of $H_2$. 
We remark that this way of expressing points in $Q_2(n)$ is \textbf{not unique}, 
and different choices of ${R}_j$ from $H_2$ may lead to the same points in 
$Q_2(n)$. 

Then one has based on (\ref{rootpoints}): 
\begin{prop}\label{B}
\begin{equation}
Q_2(n)=\left\lbrace \sum_{j=0}^{9} n_j\xi^j \mid
n_j\in \N^0, \sum_{j=0}^{9} n_j=l\leq n \right\rbrace\,.
\end{equation}
\end{prop}

Thus, $Q_2(n)$ consists of all linear combinations of up to $n$ (not necessarily different) roots from $\Delta_2$. 

\begin{proof}
It is a consequence of the fact that $T$ in (\ref{points}) is a translation by a root from $\Delta_2$ and the fact that ${R_j}$, $(j=1,\ldots,k)$, act as linear transformations on $\Delta_2$. In particular, for any tuple $(k_0,\ldots,k_9)$,   there exists a tuple $(m_0,\ldots,m_9)$ with $\sum_j k_j = \sum_j m_j$ such that  
\begin{equation}
{R_j} \sum_{j=0}^9 k_j\xi^j = \sum_{j=0}^9 m_j \xi^j\,.
\end{equation}
Thus, there exists a tuple $(n_0,\ldots,n_9)$ with 
\begin{equation}\label{exppoints}
{R}_l T {R}_{l-1} T \ldots T {R}_1 T  0 =
\sum_{j=0}^9 n_j \xi^j \,
\end{equation}
and the claim follows from (\ref{points}) since the ${R}_j$ in (\ref{exppoints}) may represent any element of $H_2$. 
\end{proof}

Hence, according to (\ref{Ldef}) and Proposition \ref{B}, we know that $L_{\alpha_H}(n)$ is the point set which corresponds to all points which are obtained by a linear combinations of up to $n$ elements from $\Delta_2$ and lie on the line $\R \alpha_H$:
\begin{equation}\label{gamma}
\begin{array}{rcl}
L_{\alpha_H}(n) & = & \R \alpha_H \cap \left\lbrace \sum_{j=0}^{9} n_j\xi^j \mid
n_j\in \N^0, \sum_{j=0}^{9} n_j=l\leq n \right\rbrace\\ 
& = & 
\lbrace \gamma \alpha_H \mid \exists (n_0,\ldots,n_9) \mbox{ such that }  
\gamma \alpha_H = \sum_{j=0}^9 n_j \xi^j, \gamma\in\R \rbrace
\end{array}
\end{equation}
The point set $L_{\alpha_H}(n)$ is thus characterized by the values $\gamma\in\R$ in (\ref{gamma}). It is 
our aim to determine the values of $\gamma$ in the following, and we aim at finding for given $n\in \N$ all $\gamma\in\R$ with  
\begin{equation}\label{gamma2}
\gamma =  \sum_{j=0}^9 n_j \xi^j 
\end{equation}
where $l\leq n$, $n_j\in \N^0$ and $\sum_{j=0}^9 n_j =l$.  

We remark that in order to facilitate notation, we will consider in the following the set $L_{\alpha_1}(n)$ where $\alpha_1=\xi^0=1$ is one of the simple roots. As mentioned before, it coincides with $L_{\alpha_H}(n)$ when viewed as a one-dimensional point set without orientation in $\R^2$, and the advantage of considering this set lies in the fact that all points are multiplied by $\alpha_1=1$ instead of $\alpha_H$. 

We start by  setting up some terminology: 
\begin{de}
Let ${\hat L_{\alpha_1}}(n):=L_{\alpha_1}(n)\setminus L_{\alpha_1}(n-1)$. 
Then we call the parameter $n$ in ${\hat L_{\alpha_1}}(n)$ the 
(growth-) level of ${L_{\alpha_1}}(n)$ and 
the points in ${\hat L_{\alpha_1}}(n)$ are called {\it points of level $n$}. 
\end{de}

Observe that the $n$th level consists of all points which are linear
combinations of exactly 
$n$ elements from $\Delta_2$, i.e. 
\begin{equation}
x\in {\hat L_{\alpha_1}}(n) 
\Leftrightarrow
x=\sum_{j=0}^9 n_j \xi^j \mbox{ with } n_j \in \N^0, \sum_{j=0}^9 n_j = n\,.
\end{equation}

Then we have 
\begin{prop}\label{C}
${\hat L_{\alpha_1}}(2)$ consists of the points $\lbrace \pm 2, \pm \tau,
\pm \tau' \rbrace$. 
\end{prop}

\begin{proof}
$\pm 2$ corresponds to $\pm 2\xi^0$, $\tau$ corresponds to $\xi^1+\xi^9$,
$-\tau$ to $\xi^4+\xi^6$, $\tau'$ to $\xi^2+\xi^8$ and $-\tau'$ to
$\xi^3+\xi^7$. No other combinations are possible. 
\end{proof}

\begin{de}
A combination $\sum_{j=0}^{9} n_j \xi^j=\gamma \in\R$ is called {\it
reducible} if it can be 
decomposed as $\gamma=\gamma_1 +\gamma_2 \in\R$ where 
$\gamma_s=\sum_{j=0}^{9} n^s_j\xi^j \in\R$, $s=1,2$ and
$\sum_{j=0}^9 (n^1_j+ n^2_j)=\sum_{j=0}^9 n_j$. Otherwise, it is called {\it nontrivial}.  
\end{de}

\begin{lem}\label{comb}
If a nontrivial combination exists on level $k\geq 3$, then it is a
combination of elements from 
$\lbrace \xi^1, \xi^4, \xi^7, \xi^8 \rbrace $. 
\end{lem}

\begin{proof}
Any combination which contains any of the pairs  
$\{\xi^1,\xi^9\}$, $\{\xi^2,\xi^8\}$, $\{\xi^3,\xi^7\}$, $\{\xi^4,\xi^6\}$,  
or, at least one of the roots $\xi^0$ and $\xi^5$, is necessarily
reducible by Proposition \ref{C} to configurations on level 2 and 1. 
Furthermore, any combination containing simultaneously $\xi^j$ and
$\xi^{(j+5) \mbox{mod} 10}$ 
is reducible to a combination on level $n-2$. 
Thus only combinations from $\lbrace \xi^1, \xi^4, \xi^7, \xi^8 \rbrace $
or 
$\lbrace \xi^2, \xi^3, \xi^6, \xi^9 \rbrace $ are potentially leading to
nontrivial combinations on level $n\geq 3$. 
Since both sets give rise to the same one-dimensional point set, 
the claim is proven. 
\end{proof}

\begin{thm}\label{main}
There is no nontrivial combination on level $n\geq 3$. 
\end{thm}

\begin{proof}
According to Lemma \ref{comb} any nontrivial combination on level $n \geq
3$ would be of the form 
\begin{equation}
\lambda_1 \xi^1+\lambda_2 \xi^4+\lambda_3 \xi^7 +\lambda_4 \xi^8, \qquad \lambda_j\in\Z, \quad j=1,\ldots,4\,.
\end{equation}
Denote the lattice spanned by $\xi^1$ and $\xi^4$ as $X$ and the one
spanned by 
$\xi^7$ and $\xi^8$ as $Y$ (see Fig. 5, 6), i.e. define 
\begin{equation}\label{lat}
\begin{array}{rcl}
X & := & \lbrace \frac{\gamma}{2} x_a + \frac{\lambda}{2} x_b \mid \gamma
\in\N, \lambda\in 
\lbrace \gamma, \gamma-1,\ldots,-(\gamma-1),-\gamma\rbrace \rbrace\\
Y & := & \lbrace \frac{\gamma}{2} y_a + \frac{\lambda}{2} y_b \mid \gamma
\in\N, \lambda\in 
\lbrace \gamma, \gamma-1,\ldots,-(\gamma-1),-\gamma\rbrace \rbrace \,.
\end{array}
\end{equation}

\begin{figure}[ht]\label{five}
\begin{center}
\includegraphics[width=3.5cm,keepaspectratio]{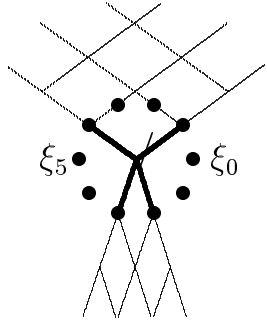}\\*
\end{center}
\caption{Displaying the lattices $X$ and $Y$.}
\end{figure}

Here $x_a$ and $x_b$, as well as $y_a$ and $y_b$ 
denote the diagonals of the parallelograms which
constitute the lattices $X$ and $Y$, respectively. 

\begin{figure}[ht]\label{six}
\begin{center}
\includegraphics[width=3.5cm,keepaspectratio]{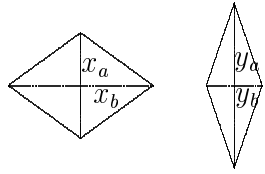}\\*
\end{center}
\caption{Displaying the defining parallelograms of the lattices $X$ and $Y$.}
\end{figure}

They serve as an orthogonal basis for the lattices $X$ and $Y$ and are of the following lengths:
\begin{equation}\label{length}
\begin{array}{rl}
l(x_a)=\sqrt{3-\tau} & l(x_b)=\tau\\
l(y_a)=\sqrt{2+\tau} & l(y_b)=\tau -1\,.
\end{array}
\end{equation}

A necessary condition for a nontrivial combination to exist is thus that there exist $\gamma_j$ and $\lambda_j$ as in
(\ref{lat}) such that  
\begin{equation}\label{cond}
\begin{array}{rcl}
\frac{\gamma_1}{2}l(x_a) & = & \frac{\gamma_2}{2} l(x_b)\\
\frac{\lambda_1}{2}l(y_a) & = & \frac{\lambda_2}{2} l(y_b)\,.
\end{array}
\end{equation}
However, this implies $\lambda_1=\lambda_2=0$ and $\gamma_1=\gamma_2=0$,  
which proves the claim. 
\end{proof}

Based on Theorem \ref{main} we have 
\begin{coro}
\begin{equation}\label{PointsB}
L_{\alpha_1}(n)= \lbrace (a+c)+(b-c)\tau \mid a,b,c\in\Z, |a|+2|b|+2|c|\leq n
\rbrace
\end{equation}
\end{coro}

\begin{proof}
Follows via $a+b\tau+c\tau'= (a+c)+(b-c)\tau$ from a decomposition of each
level into the contributions from level 2 and 1. 
\end{proof}

Note, that correspondingly 
\begin{equation}\label{PointsC}
L_{\alpha_H}(n)= \lbrace ((a+c)+(b-c)\tau) \alpha_H \mid a,b,c\in\Z, |a|+2|b|+2|c|\leq n
\rbrace
\end{equation}
describes the dominant points $(a b)^T$ with $a=b$. 

\subsubsection{Comparison with cut-and-project quasicrystals:}

The advantage of expressing $L_{\alpha_1}(n)$ as in (\ref{PointsB}) is the fact that it facilitates comparison with the cut-and-project scheme. 
For this purpose, we briefly recall the definition of a one-dimensional cut-and-project quasicrystal associated with the irrationality $\tau$ ( see \cite{Patera:1999} and references within): 

Consider the algebraic number field $\mathbbm{Q}[\sqrt{5}]$\ and 
its nontrivial automorphism denoted by~$'$ and defined by 
$a+b\sqrt{5}\rightarrow a-b\sqrt{5}$ with $a,b\in\Z$. 
In particular, ~$'$ transforms $\tau$ into 
$\tau'=\frac12 (1-\sqrt{5})$. Furthermore, denote the ring of 
integers of $\Q[\sqrt{5}]$ by $\Z[\tau ]= \Z+\Z\tau$. 
Then we have: 
\begin{de}\label{cutpro}
Let $\Omega$ be a bounded interval. The point set 
\begin{equation} 
\label{quasic} 
\Sigma(\Omega):= \left\{ x\in\Z[\tau]\mid\, x'\in \Omega 
\right\}\, 
\end{equation} 
is called cut-and-project quasicrystal, and 
the interval $\Omega$ is called the 
acceptance window of $\Sigma(\Omega)$. 
\end{de} 
Based on this, we obtain: 
\begin{prop}\label{relA}
\begin{equation}\label{rel1}
L_{\alpha_1}(n)\subset \Sigma([-n,n])\cap [-n,n]
\end{equation}
\end{prop}

\begin{proof}
Clearly, $x\in L_{\alpha_1}(n)$ implies $x\in\Z[\tau]$. Furthermore, 
$|x'|\leq n$ and $|x|\leq n$, thus $x\in \Sigma([-n,n])\cap [-n,n]$.
\end{proof}

Note that the opposite inclusion does {\it not} hold, so that the two sets 
are not equal: 

\begin{lem}
The inclusion in Proposition \ref{relA} is a true inclusion. Deficiencies
occur for $n\geq 3$. 
\end{lem}

\begin{proof}
Let $x:=x_1+\tau x_2\in \Sigma([-n,n])\cap [-n,n]$ and suppose
w.l.o.g. that $x_2>0$. 
Then $x_1$ is bounded by 
\begin{equation}\label{bound} 
-n-\tau' x_2\leq x_1 \leq n-\tau x_2\,.
\end{equation}
Fix the $x_2$-component. Then a sufficient condition for the existence of an $x_1\in\integer$ fulfilling (\ref{bound})
is 
\begin{equation}
n-\tau x_2 -1 \geq -n-\tau' x_2
\end{equation}
which implies 
\begin{equation}
x_2\leq \left\lbrack \frac{2n-1}{\tau-\tau'}\right\rbrack =: N_n\,.
\end{equation}
On the other hand, for $L_n:=\lbrace (b,c)\in\Z\times\Z \mid 2|b|+2|c|\leq n \rbrace$ we
have  
\begin{equation}
\max_{(b,c)\in L_n} (b-c)=: M_n
\end{equation}
where 
\begin{equation}
M_n=\left\lbrace 
\begin{array}{ll}
\frac{n}{2} & n \mbox{ even }\\
\frac{n-1}{2} & n \mbox{ odd }\,.
\end{array}
\right.
\end{equation}
Since $M_n < N_n$ for $n\geq 3$, deficiencies occur. 
\end{proof}

\noindent
We remark that for $n=1,2$, the point sets coincide and that deficiencies
indeed occur only for $n \geq 3$. 

\noindent
\textbf{Example:}

\noindent
In the case  $n=3$, $N_3=2$ and $M_3=1$, which is consistent with $-1+2\tau \in
\Sigma([-3,3])\cap [-3,3]$ 
but $-1+2\tau\not\in L_{\alpha_1}(3)$.

\begin{coro}
$L_{\alpha_1}(n)$ does not correspond to a cut-and-project quasicrystal
with connected acceptance window for $n\geq 3$. 
\end{coro}

Note that as a consequence of Proposition \ref{relA} we obtain a lower bound for the minimal distance 
between adjacent points in $L_{\alpha_1}(n)$: 

\begin{lem}\label{Min}
The minimal distance in $L_{\alpha_1}(n)$ is greater or equal to the one in $\Sigma([-n,n])$.  
\end{lem}

We remark that the latter has been determined in \cite{Patera:1998b} and varies in dependence on the size of the acceptance window.

Finally, let us make a remark about the {\it{repetitivity}} and {\it
scaling} properties of patterns $P$ in $L_{\alpha_1}(n)$. 
As follows from (\ref{PointsB}) we have:  

\begin{itemize}
\item
For any pattern $P$ with $P\subset L_{\alpha_1}(r)$ and $x\in L_{\alpha_1}(s)$
we have 
$(P+x)\subset L_{\alpha_1}(r+s)$. Thus, multiple pattern repetitions occur
with growing $n$. 
\item 
For any $n\in\N$ there exists $l\in\N$, $l\geq 2n$, such that $\tau L_{\alpha_1}(n) \subset
L_{\alpha_1}(l)$, as follows from  
\begin{equation}
\tau L_{\alpha_1}(n)= \lbrace (b-c)+(a+b)\tau \mid a,b,c\in\Z,
|a|+2|b|+2|c|\leq n \rbrace\,.
\end{equation}
\end{itemize}

\subsection{Implications for dominant points $(a\ b)^T$ with $a\not= b$}

In this subsection, we use the information on the dominant points $(a\ b)^T$ with $a= b$ 
derived previously in order to infer information also 
about the $a\not= b$ case. 

We start by showing: 
\begin{thm}
For each $x\in Q_2(n)$ there exist $y\in L_{\xi^0}(n)$ and $z\in
L_{\xi^1}(n)$ such that $x=y+z$. 
\end{thm}

Note that the above statement is trivial if we replace $Q_2(n)$ by
$Q_2(2n)$ and needs to be proven only for the case that the cut-off
parameter of the ($L$-)subspaces coincides with the cut-off of the
two-dimensional ($Q$-)setting. 

\begin{proof}
Let $x\in Q_2(n)$. Then $x=\sum_{j=0}^{4} \beta_j \xi^j$, where $\sum_{j=0}^{4}
|\beta_j| \leq n$ with $\beta_j\in\Z$ and $\xi^j\in\Delta_2$ as in (\ref{rootpoints}). 
Expressing $\xi^j$ for $j\geq 2$ in terms of $\xi^0$ and $\xi^1$ leads to 
\begin{equation}\label{x}
x=\lbrace (\beta_0 - \beta_2)-\tau (\beta_3 + \beta_4)\rbrace \xi^0
+\lbrace (\beta_1 +\beta_4)+\tau (\beta_2 + \beta_3) \rbrace \xi^1\,.
\end{equation}
On the other hand, 
\begin{equation}
\label{yz}
y+z = \lbrace (a_1+c_1)+\tau (b_1 -c_1)\rbrace \xi^0
+\lbrace (a_2+c_2)+\tau (b_2 -c_2) \rbrace \xi^1
\end{equation}
where $|a_j|+2 |b_j|+2|c_j| \leq n$. 

We thus have to show that for all $\beta_j$ with $\sum_{j=0}^{4} |\beta_j| \leq n$
there exist $a_j$, $b_j$ and $c_j$ with 
$|a_j|+2 |b_j|+2 |c_j|\leq n$ such that the following equalities hold: 

\begin{equation}
\label{equals}
\begin{array}{rclrcl}
\beta_0 - \beta_2 & = & a_1 + c_1 & -(\beta_3 + \beta_4) & = & b_1 - c_1
\\
\beta_1+\beta_4 & = & a_2 + c_2 & \beta_2 + \beta_3 & = & b_2 - c_2 \,.
\end{array}
\end{equation}

With the definitions 
\begin{equation}
\label{fg}
\begin{array}{rcl}
f(c_1) & := & |a_1(\beta_0,\beta_2,c_1)|+2 |b_1(\beta_3,\beta_4,c_1)|+2 |c_1| \\
& = & |\beta_0 -\beta_2-c_1| + 2 |c_1-(\beta_3+\beta_4)|+2|c_1|\\
g(c_2) & := & |a_2(\beta_1,\beta_4,c_2)|+2 |b_2(\beta_2,\beta_3,c_2)|+2 |c_2|\\
& = & |\beta_1 +\beta_4-c_2| + 2 |c_2+\beta_2+\beta_3|+2|c_2|
\end{array}
\end{equation}
this is equivalent to showing that for all $\beta_j$ with $\sum_{j=0}^{4} |\beta_j|
\leq n$ there exist $c_1$ and $c_2$ such that $f(c_1)\leq n$ and
$g(c_2)\leq n$. 

For this, we investigate minima and maxima of these functions in
dependence on the parameter ranges. In particular, we have 
\begin{equation}
\label{fgprime}
\begin{array}{rcl}
f'(c_1)  &= & 
\left\lbrace 
\begin{array}{rl}
-1 & \beta_0-\beta_2 > c_1\\
1  & \beta_0-\beta_2 < c_1
\end{array}
\right\rbrace
+ 2
\left\lbrace 
\begin{array}{rl}
-1 & \beta_3+\beta_4 > c_1\\
1  & \beta_3+\beta_4 < c_1
\end{array}
\right\rbrace
+ 2
\left\lbrace 
\begin{array}{rl}
-1 & 0 > c_1\\
1  & 0 < c_1
\end{array}
\right\rbrace\\
g'(c_2)  &= & 
\left\lbrace 
\begin{array}{rl}
-1 & \beta_1+\beta_4 > c_2\\
1  & \beta_1+\beta_4 < c_2
\end{array}
\right\rbrace
+ 2
\left\lbrace 
\begin{array}{rl}
-1 & -(\beta_2+\beta_3) > c_2\\
1  & -(\beta_2+\beta_3) < c_2
\end{array}
\right\rbrace
+ 2
\left\lbrace 
\begin{array}{rl}
-1 & 0 > c_2\\
1  & 0 < c_2
\end{array}
\right\rbrace\,.
\end{array}
\end{equation}
The choices of parameters leading to different qualitative behaviour
of $f'(c_1)$ and $g'(c_2)$ are discussed separately. 
For instance, for $0\leq \beta_3+\beta_4 < \beta_0 - \beta_2$ the minimum
of the function $f(c_1)$ is at $c_1=\beta_3+\beta_4$ and we have 
\begin{equation}
\begin{array}{rcl}
f(\beta_3+\beta_4) & = &
|\beta_0-\beta_2-\beta_3-\beta_4|+2|\beta_3+\beta_4|\\
&= & \beta_0-\beta_2+\beta_3+\beta_4 \leq \sum |\beta_j| \leq n\,.
\end{array}
\end{equation}
The other cases can be treated analogously, which proves the claim. 
\end{proof}

Let $C(\xi^0,\xi^1)$ denote the cone enclosed by the halflines $\R^+ \xi^0$ and $\R^+ \xi^1$. Then this theorem shows that any point in $Q_2(n) \cap C(\xi^0,\xi^1)$ can be expressed as a linear combination of points from $L_{\xi^0}(n)$ and $L_{\xi^1}(n)$ when viewed as vectors in $\R^2$. Since the points in $Q_2(n) \cap C(\xi^0,\xi^1)$ describe the whole point set $Q_2(n)$ due to 10-fold rotational symmetry, it follows that any dominant point in $Q_2(n)$ can be expressed as a linear combination of points from $L_{\xi^j}(n)$ and $L_{\xi^{j+1}}(n)$ for suitably chosen $j\in\{0,\ldots,9\}$. In particular, dominant points $(a b)^T$ with $a>b$ are given by linear combinations from $L_{\xi^1}(n)$ and $L_{\xi^2}(n)$, and 
dominant points $(a b)^T$ with $a<b$ are given by linear combinations from $L_{\xi^2}(n)$ and $L_{\xi^3}(n)$.

We remark that some of the properties proven here are rooted in the special structure of the ring of cyclotomic integers, which based on (\ref{rootpoints}) is given by 
\begin{equation}
\integer[\xi] = \sum_{j=0}^9 \integer \xi^j = \integer [\tau] + \integer [\tau] \xi
\end{equation}
and of which $Q_2(n)$ is by construction a subset. 

We can again embed our point set into a cut-and-project quasicrystal. 
For this, we indicate briefly how the setting of cut-and-project quasicrystals as introduced in Definition \ref{cutpro} can be generalized to two dimensional point sets with $H_2$ symmetry (see \cite{P} and references within for more details): 

\begin{de}\label{star}
Let $M$ denote a $\zt$-lattice with respect to some basis in $\R^k$. Then we call the map 
$^*:M\mapsto \R^k$ with the property  
$(ax+y)^*=a'x^*+y^*$ for all $x$, $y\in M$ and $a\in\zt$ a $*$-map. 
\end{de}

\begin{de}\label{cutpro2D}
Let $M$ be a $\zt$-lattice in $\R^k$ and $\Omega$ a bounded region in $\R^k$, called acceptance window. 
Then 
\begin{equation}
\Sigma(\Omega)=\left\lbrace x\in M|x^*\in\Omega\right\rbrace
\end{equation}
defines a cut-and-project quasicrystal in $k$ dimensions.
\end{de}
 
In order to construct a cut-and-project quasicrystal with $H_2$ symmetry along these lines, 
one takes $M={\zt} \Delta_2$ as $\zt$-lattice in $\R^k$ in Def. \ref{cutpro2D}, and as 
${}^*$-map, one uses based on (\ref{rootpoints}) 
\begin{equation}
{}^*: \xi \mapsto \xi^2 \,,
\end{equation} 
which fulfills the requirements of a $*$-map (cf. Def. \ref{star}) and leaves $\Delta_2$ invariant:  
\begin{equation}
{}^*:\Delta_2\rightarrow\Delta_2^*\equiv \Delta_2\,.
\end{equation}
With this, a cut-and-project quasicrystal as in (\ref{twoQC}) can be parameterized as 
\begin{equation}\label{twoQC}
\Sigma(\Omega)  =  \left\lbrace (x_1+\tau x_2)\alpha_1+(x_3+\tau x_4)\alpha_2 \right.
\mid \left.(x_1+\tau' x_2)\alpha_1^* +(x_3+\tau' x_4)\alpha_2^* \in\Omega\right\rbrace\, 
\end{equation}
where $\Omega$ may be any bounded region in $\R^2$. 

Note in particular that $\alpha_1^*=(\xi^0)^*=\xi^0$ and $\alpha_2^*=(\xi^4)^*=\xi^8$ (compare with (\ref{rootpoints})). 

Based on this, we find the following in our context:
\begin{lem}\label{cut}
Let $D(n)$ denote the convex hull of the regular decagon inscribed into a circle of radius $n$ around the origin and let $\Sigma(D(n))$ denote the corresponding cut-and-project quasicrystal. 
Then  
\begin{equation}
Q_2(n) \subset \Sigma (D(n)) \cap D(n)\,.
\end{equation}
\end{lem}

\begin{proof}
$x\in Q_2(n)$ implies $x=\sum_{j=0}^{9} n_j \xi^j$ with $\xi^j\in
\Delta_2$, $n_j\in \N^0$ and $\sum_{j=0}^{9} n_j = l \leq n$. Thus, $x\in D(n)$. Since $*$ leaves $\Delta_2$ invariant, $x^*\in D(n)$ and the
claim follows. 
\end{proof}

As before, there is no equality in Lemma \ref{cut} and for $n\geq 3$ deficiencies occur. Also, minimal distances in $Q_2(n)$ are bounded from below by the ones in $\Sigma (D(n))$. 

\section{Conclusion}

We have suggested a new way to construct mathematical models for fragments of aperiodic 
point sets with 10-fold symmetry. Like cut-and-project
quasicrystals, they require some cut-off condition which prevents the sets 
from becoming dense. The special feature of these models is that they are
by construction finite structures -- not idealized infinite ones -- which
grow from a seed point as we demonstrate in Appendix B. We have shown that they are not fragments of sets obtainable via the cut-and-project scheme for convex windows. The restriction to convex windows is a plausible assumption when dealing with growth processes which are not hindered by obstacles. In the presence of obstacles, points would be generated around them and the corresponding acceptance windows would not necessarily be convex and connected. The deviation of our point sets from the cut-and-project situation with convex windows is given in terms of  a set of ``deficiencies'', which appear close to the boundary of the set. The occurrence of deficiencies is a novel aspect and a special feature of our models, and it has to be discussed in how far it may help to model growth deficiencies which occur in real life quasicrystals. 

Finally, we remark that though our initial motivation for this study comes from the field of quasicrystals, we expect that the mathematical structures provided by the affine extension of noncrystallographic Coxeter groups will open also other fields of applications. We plan to investigate in a next step the application of our results to the study of fullerenes, in particular the description of onion like structures and nanotubes.

\subsection*{Acknowledgements} 
\par
J.~P. acknowledges financial support by the Natural Sciences and
Engineering 
 Research Council of Canada and FCAR of Quebec and R.~T. by a Marie Curie fellowship.  
She is grateful for the hospitality extended to her 
at the Centre de Recherches Math\'ematiques, Universit\'e de 
Montr\'eal, where this work has been started.



\section*{Appendix A}

In this appendix, we indicate the generalized Cartan matrices obtained
after relaxation of the condition $(a_{ij})\in\integer[\tau]^-$ (compare
with Section 3). 
In particular, let $a:=a_1+\tau a_2$, $b:=b_1+\tau b_2$, $c:=c_1+\tau c_2$
and $d:=d_1+\tau d_2$ be the entries of the matrices 
\begin{equation}\label{genCarts}
\left( \begin{array}{ccc}
 2 & a & b \\ 
 a & 2 & -\tau \\ 
 b & -\tau & 2
\end{array}\right)
,
\left( \begin{array}{cccc}
 2    &  a  &   b   & c     \cr 
 a    &  2  &  -1   & 0     \cr
 b    &  -1 &   2   & -\tau \cr
 c    &  0  & -\tau &   2  
\end{array}\right)
\mbox{ and } 
\left( \begin{array}{ccccc}
 2    &    a   &  b  &   c   &  d      \cr 
 a    &    2   &  -1 &   0   &  0      \cr
 b    &   -1   &   2 &   -1  &  0      \cr
 c    &    0   &  -1 &    2  & -\tau   \cr 
 d    &    0   &  0  & -\tau &  2
\end{array}\right)\,.
\end{equation}
Then the entries in the following tables define generalized Cartan
matrices 
for $H_2$ (first matrix in (\ref{genCarts}) and Table 1), $H_3$ (second
matrix in (\ref{genCarts}) and Table 2) and $H_4$ (third matrix in
(\ref{genCarts}) and Table 3), respectively.

\textbf{Table 1: the case of $H_2$}
\begin{longtable}{|rr|rr|}
\hline$a_1$&$a_2$&$b_1$&$b_2$\\\hline\endhead
\hline\endfoot
-2 & 0 & 0 & 1\\
-1 & 1 & -1 & 1\\
-1 & 1 & 0 & -1\\
0 & -1 & -1 & 1\\
0 & -1 & 2 & 0\\
0 & 1 & -2 & 0\\
0 & 1 & 1 & -1\\
1 & -1 & 0 & 1\\
1 & -1 & 1 & -1\\
2 & 0 & 0 & -1\\
\end{longtable}

\textbf{Table 2: the case of $H_3$}
\begin{longtable}{|rr|rr|rr|}
\hline$a_1$&$a_2$&$b_1$&$b_2$&$c_1$&$c_2$\\\hline\endhead
\hline\endfoot
-2 & 0 & 1 & 0 & 0 & 0\\
-1 & -2 & 0 & 2 & -2 & 1\\
-1 & 1 & 0 & -1 & 1 & 0\\
-1 & 1 & 0 & 0 & -1 & 0\\
-1 & 1 & 1 & -1 & 1 & 0\\
-1 & 1 & 1 & 0 & -1 & 0\\
0 & -3 & -1 & 2 & -1 & 1\\
0 & -1 & -1 & 1 & -1 & 1\\
0 & 0 & -1 & 1 & 0 & 0\\
0 & 0 & 0 & -1 & 2 & 0\\
0 & 0 & 0 & 1 & -2 & 0\\
0 & 0 & 1 & -1 & 0 & 0\\
0 & 1 & 1 & -1 & 1 & -1\\
0 & 3 & 1 & -2 & 1 & -1\\
1 & -1 & -1 & 0 & 1 & 0\\
1 & -1 & -1 & 1 & -1 & 0\\
1 & -1 & 0 & 0 & 1 & 0\\
1 & -1 & 0 & 1 & -1 & 0\\
1 & 2 & 0 & -2 & 2 & -1\\
2 & 0 & -1 & 0 & 0 & 0\\
\end{longtable}

\textbf{Table 3: the case of $H_4$}
\begin{longtable}{|rr|rr|rr|rr|}
\hline
$a_1$&$a_2$&$b_1$&$b_2$&$c_1$&$c_2$&$ d_1$&$d_2$\\
\hline\endhead
\hline\endfoot
-2 & 0 & 1 & 0 & 0 & 0 & 0 & 0\\
-1 & 0 & -1 & 0 & 1 & 0 & 0 & 0\\
-1 & 0 & 0 & 0 & -1 & 0 & 0 & 1\\
-1 & 0 & 0 & 0 & 0 & 0 & -1 & 1\\
-1 & 0 & 0 & 0 & 0 & 1 & 0 & -1\\
-1 & 0 & 0 & 0 & 1 & -1 & 1 & 0\\
-1 & 0 & 0 & 0 & 1 & 0 & -1 & 0\\
-1 & 0 & 0 & 0 & 1 & 0 & 1 & -1\\
-1 & 0 & 0 & 1 & 0 & -1 & 1 & 0\\
-1 & 0 & 0 & 1 & 0 & 0 & -1 & 0\\
-1 & 0 & 0 & 1 & 1 & -1 & 0 & 0\\
-1 & 0 & 1 & -1 & -1 & 1 & 0 & 0\\
-1 & 0 & 1 & -1 & 0 & 0 & 1 & 0\\
-1 & 0 & 1 & -1 & 0 & 1 & -1 & 0\\
-1 & 0 & 1 & 0 & -1 & 0 & -1 & 1\\
-1 & 0 & 1 & 0 & -1 & 0 & 1 & 0\\
-1 & 0 & 1 & 0 & -1 & 1 & -1 & 0\\
-1 & 0 & 1 & 0 & 0 & -1 & 0 & 1\\
-1 & 0 & 1 & 0 & 0 & 0 & 1 & -1\\
-1 & 0 & 1 & 0 & 1 & 0 & 0 & -1\\
-1 & 0 & 2 & 0 & -1 & 0 & 0 & 0\\
-1 & 1 & 0 & -1 & 0 & 0 & 1 & 0\\
-1 & 1 & 0 & -1 & 0 & 1 & -1 & 0\\
-1 & 1 & 0 & -1 & 1 & 0 & 0 & 0\\
-1 & 1 & 0 & 0 & 0 & -1 & 0 & 1\\
-1 & 1 & 0 & 0 & 0 & 0 & 0 & 0\\
-1 & 1 & 0 & 0 & 1 & 0 & 0 & -1\\
-1 & 1 & 1 & -1 & -1 & 0 & 0 & 1\\
-1 & 1 & 1 & -1 & 0 & 0 & 0 & 0\\
-1 & 1 & 1 & -1 & 0 & 1 & 0 & -1\\
-1 & 1 & 1 & 0 & -1 & 0 & 0 & 0\\
-1 & 1 & 1 & 0 & 0 & -1 & 1 & 0\\
-1 & 1 & 1 & 0 & 0 & 0 & -1 & 0\\
0 & -1 & -1 & 1 & 0 & 0 & 0 & 0\\
0 & -1 & -1 & 1 & 1 & -1 & 1 & 0\\
0 & -1 & -1 & 1 & 1 & 0 & -1 & 0\\
0 & -1 & 0 & 0 & -1 & 1 & 0 & 0\\
0 & -1 & 0 & 0 & 0 & 0 & -1 & 1\\
0 & -1 & 0 & 0 & 1 & 0 & 1 & -1\\
0 & -1 & 0 & 1 & -1 & 0 & -1 & 1\\
0 & -1 & 0 & 1 & 0 & 0 & 1 & -1\\
0 & -1 & 0 & 1 & 1 & -1 & 0 & 0\\
0 & -1 & 1 & 0 & -1 & 0 & 1 & 0\\
0 & -1 & 1 & 0 & -1 & 1 & -1 & 0\\
0 & -1 & 1 & 0 & 0 & 0 & 0 & 0\\
0 & 0 & -1 & 0 & 0 & 0 & -1 & 1\\
0 & 0 & -1 & 0 & 0 & 0 & 1 & 0\\
0 & 0 & -1 & 0 & 0 & 1 & -1 & 0\\
0 & 0 & -1 & 0 & 1 & -1 & 0 & 1\\
0 & 0 & -1 & 0 & 1 & 0 & 1 & -1\\
0 & 0 & -1 & 0 & 2 & 0 & 0 & -1\\
0 & 0 & -1 & 1 & 0 & -1 & 0 & 1\\
0 & 0 & -1 & 1 & 1 & -1 & 0 & 0\\
0 & 0 & -1 & 1 & 1 & 0 & 0 & -1\\
0 & 0 & 0 & -1 & -1 & 1 & -1 & 1\\
0 & 0 & 0 & -1 & 0 & 1 & 1 & -1\\
0 & 0 & 0 & -1 & 1 & 0 & 0 & 0\\
0 & 0 & 0 & 0 & -1 & 0 & 1 & 0\\
0 & 0 & 0 & 0 & -1 & 1 & -1 & 0\\
0 & 0 & 0 & 0 & 0 & -1 & 2 & 0\\
0 & 0 & 0 & 0 & 0 & 1 & -2 & 0\\
0 & 0 & 0 & 0 & 1 & -1 & 1 & 0\\
0 & 0 & 0 & 0 & 1 & 0 & -1 & 0\\
0 & 0 & 0 & 1 & -1 & 0 & 0 & 0\\
0 & 0 & 0 & 1 & 0 & -1 & -1 & 1\\
0 & 0 & 0 & 1 & 1 & -1 & 1 & -1\\
0 & 0 & 1 & -1 & -1 & 0 & 0 & 1\\
0 & 0 & 1 & -1 & -1 & 1 & 0 & 0\\
0 & 0 & 1 & -1 & 0 & 1 & 0 & -1\\
0 & 0 & 1 & 0 & -2 & 0 & 0 & 1\\
0 & 0 & 1 & 0 & -1 & 0 & -1 & 1\\
0 & 0 & 1 & 0 & -1 & 1 & 0 & -1\\
0 & 0 & 1 & 0 & 0 & -1 & 1 & 0\\
0 & 0 & 1 & 0 & 0 & 0 & -1 & 0\\
0 & 0 & 1 & 0 & 0 & 0 & 1 & -1\\
0 & 1 & -1 & 0 & 0 & 0 & 0 & 0\\
0 & 1 & -1 & 0 & 1 & -1 & 1 & 0\\
0 & 1 & -1 & 0 & 1 & 0 & -1 & 0\\
0 & 1 & 0 & -1 & -1 & 1 & 0 & 0\\
0 & 1 & 0 & -1 & 0 & 0 & -1 & 1\\
0 & 1 & 0 & -1 & 1 & 0 & 1 & -1\\
0 & 1 & 0 & 0 & -1 & 0 & -1 & 1\\
0 & 1 & 0 & 0 & 0 & 0 & 1 & -1\\
0 & 1 & 0 & 0 & 1 & -1 & 0 & 0\\
0 & 1 & 1 & -1 & -1 & 0 & 1 & 0\\
0 & 1 & 1 & -1 & -1 & 1 & -1 & 0\\
0 & 1 & 1 & -1 & 0 & 0 & 0 & 0\\
1 & -1 & -1 & 0 & 0 & 0 & 1 & 0\\
1 & -1 & -1 & 0 & 0 & 1 & -1 & 0\\
1 & -1 & -1 & 0 & 1 & 0 & 0 & 0\\
1 & -1 & -1 & 1 & 0 & -1 & 0 & 1\\
1 & -1 & -1 & 1 & 0 & 0 & 0 & 0\\
1 & -1 & -1 & 1 & 1 & 0 & 0 & -1\\
1 & -1 & 0 & 0 & -1 & 0 & 0 & 1\\
1 & -1 & 0 & 0 & 0 & 0 & 0 & 0\\
1 & -1 & 0 & 0 & 0 & 1 & 0 & -1\\
1 & -1 & 0 & 1 & -1 & 0 & 0 & 0\\
1 & -1 & 0 & 1 & 0 & -1 & 1 & 0\\
1 & -1 & 0 & 1 & 0 & 0 & -1 & 0\\
1 & 0 & -2 & 0 & 1 & 0 & 0 & 0\\
1 & 0 & -1 & 0 & -1 & 0 & 0 & 1\\
1 & 0 & -1 & 0 & 0 & 0 & -1 & 1\\
1 & 0 & -1 & 0 & 0 & 1 & 0 & -1\\
1 & 0 & -1 & 0 & 1 & -1 & 1 & 0\\
1 & 0 & -1 & 0 & 1 & 0 & -1 & 0\\
1 & 0 & -1 & 0 & 1 & 0 & 1 & -1\\
1 & 0 & -1 & 1 & 0 & -1 & 1 & 0\\
1 & 0 & -1 & 1 & 0 & 0 & -1 & 0\\
1 & 0 & -1 & 1 & 1 & -1 & 0 & 0\\
1 & 0 & 0 & -1 & -1 & 1 & 0 & 0\\
1 & 0 & 0 & -1 & 0 & 0 & 1 & 0\\
1 & 0 & 0 & -1 & 0 & 1 & -1 & 0\\
1 & 0 & 0 & 0 & -1 & 0 & -1 & 1\\
1 & 0 & 0 & 0 & -1 & 0 & 1 & 0\\
1 & 0 & 0 & 0 & -1 & 1 & -1 & 0\\
1 & 0 & 0 & 0 & 0 & -1 & 0 & 1\\
1 & 0 & 0 & 0 & 0 & 0 & 1 & -1\\
1 & 0 & 0 & 0 & 1 & 0 & 0 & -1\\
1 & 0 & 1 & 0 & -1 & 0 & 0 & 0\\
2 & 0 & -1 & 0 & 0 & 0 & 0 & 0\\
\end{longtable}


\section*{Appendix B}

In this appendix we demonstrate the growth of $Q_2(n)$ in dependence on $n$ for
$n=1,\ldots,6$. 
Note that we display the point sets in a circle of a radius corresponding
to four times the root length. Thus, the complete point set is visible
only until iteration step $n=4$ and is truncated afterwards. Note that since the 
point set $Q_2(2)$ is displayed in Fig. 4 we omit it in this list.

\begin{figure}[ht]\label{seven}
\begin{center}
\includegraphics[width=5.0cm,keepaspectratio]{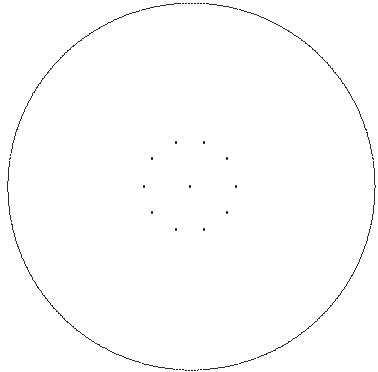}\\*
\end{center}
\caption{The point set $Q_2(1)$.}
\end{figure}

\begin{figure}[ht]\label{eight}
\begin{center}
\includegraphics[width=5.0cm,keepaspectratio]{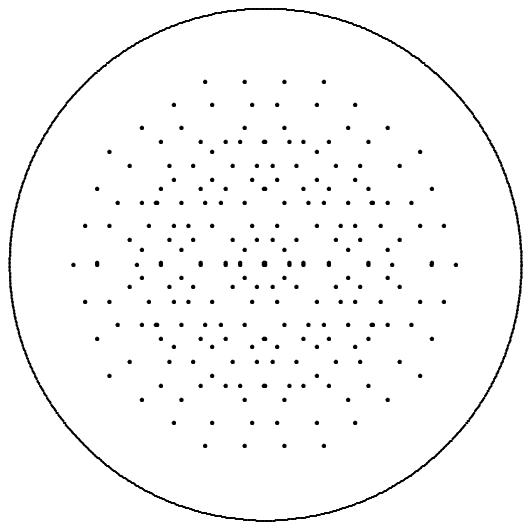}\\*
\end{center}
\caption{The point set $Q_2(3)$.}
\end{figure}

\begin{figure}[ht]\label{nine}
\begin{center}
\includegraphics[width=5.0cm,keepaspectratio]{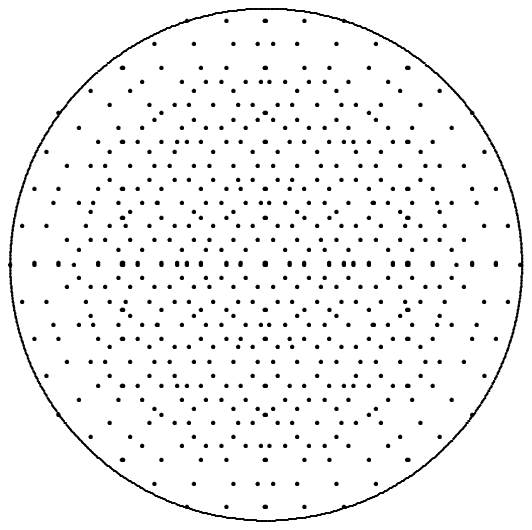}\\*
\end{center}
\caption{The point set $Q_2(4)$.}
\end{figure}

\begin{figure}[ht]\label{ten}
\begin{center}
\includegraphics[width=5.0cm,keepaspectratio]{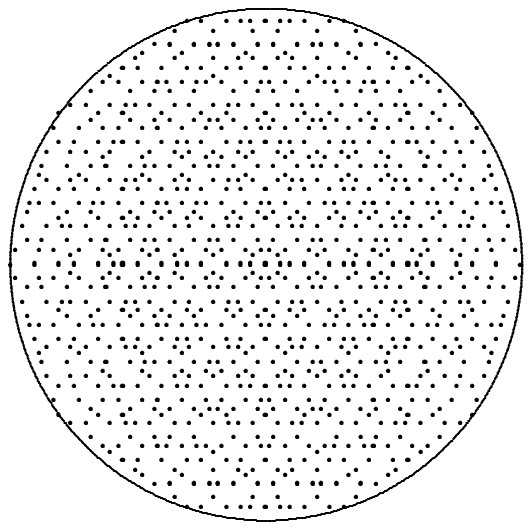}\\*
\end{center}
\caption{The point set $Q_2(5)$.}
\end{figure}

\begin{figure}[ht]\label{eleven}
\begin{center}
\includegraphics[width=5.0cm,keepaspectratio]{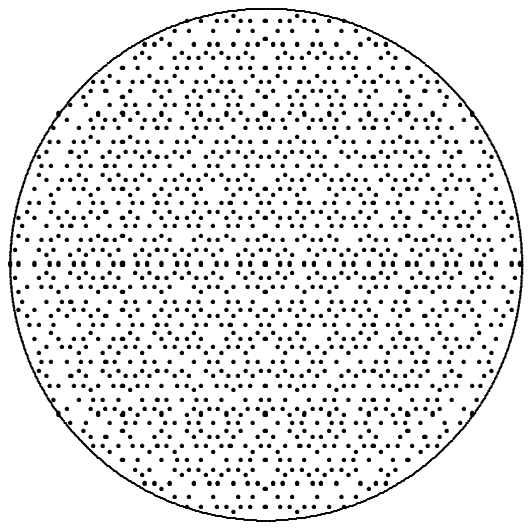}\\*
\end{center}
\caption{The point set $Q_2(6)$.}
\end{figure}

\end{document}